\documentclass{article}
\usepackage{amsthm,amsmath,amssymb,latexsym,amscd}

\renewcommand{\b}{\beta}

\newcommand{\A}{\mathcal{A}}

\newcommand{\T}{\mathcal{T}}

\newcommand{\Hom}{\mathrm{Hom}}

\newcommand{\Dend}{\mathrm{\mathbf{Dend}}}
\newcommand{\Bax}{\mathrm{\mathbf{GRB}}}
\newcommand{\phiBax}{\mathrm{\mathbf{TRB}}}
\newcommand{\NS}{\mathrm{\mathbf{NS}}}

\newcommand{\p}{\prime}
\newcommand{\ti}{\tilde}
\renewcommand{\c}{\circ}
\newcommand{\ot}{\otimes}
\newcommand{\pr}{\prec}
\newcommand{\su}{\succ}

\newtheorem{definition}{Definition}[section]
\newtheorem{lemma}[definition]{Lemma}
\newtheorem{proposition}[definition]{Proposition}

\newtheorem{corollary}[definition]{Corollary}
\newtheorem{remark}[definition]{Remark}
\newtheorem{example}[definition]{Example}
\date{}

\begin{document}
\title{
Quantum Analogy of Poisson Geometry,
Related Dendriform Algebras
and Rota-Baxter Operators}
\author{UCHINO, Kyousuke
\footnote{Post doctoral student supported
by Japan Society for the Promotion of Science}}
\maketitle

\begin{abstract}
We will introduce an associative (or quantum) version
of Poisson structure tensors.
This object is defined as an operator satisfying
a ``generalized" Rota-Baxter identity of weight zero.
Such operators are called generalized Rota-Baxter operators.
We will show that generalized Rota-Baxter operators
are characterized by a cocycle condition
so that Poisson structures are so.
By analogy with twisted Poisson structures,
we propose a new operator ``twisted Rota-Baxter operators"
which is a natural generalization of generalized Rota-Baxter operators.
It is known that classical Rota-Baxter operators
are closely related with dendriform algebras.
We will show that twisted Rota-Baxter operators
induce NS-algebras which is a twisted version of
dendriform algebra.
The twisted Poisson condition is considered
as a Maurer-Cartan equation up to homotopy.
We will show the twisted Rota-Baxter condition also is so.
And we will study a Poisson-geometric reason, how
the twisted Rota-Baxter condition arises.
\end{abstract}

\section{Introduction.}

Let $V$ be a smooth manifold. A bivector field $\pi$ on $V$
defines a Poisson bracket on the smooth functions
if and only if $\pi$ is a solution of
the classical Yang-Baxter type, or Maurer-Cartan equation,
\begin{equation}\label{cybe}
\frac{1}{2}[\pi,\pi]=0,
\end{equation}
where the bracket $[\cdot,\cdot]$ is a Schouten bracket on $V$.
$\pi$ is a solution of (\ref{cybe})
if and only if $\pi$ satisfies the operator identity (\ref{cybe2}) below.
For any smooth 1-forms $\eta,\theta$,
\begin{equation}\label{cybe2}
[\pi(\eta),\pi(\theta)]=\pi([\pi(\eta),\theta]+[\eta,\pi(\theta)])
\end{equation}
where $\pi$ is identified with
a bundle map from the cotangent bundle to the tangent bundle.
In \cite{Se}, Severa and Weinstein gave
an extended framework of Poisson geometry.
A twisted Poisson structure is to be a
solution $\pi:T^{*}V\to TV$ of the following equation,
\begin{equation}\label{tcybep}
[\pi(\eta),\pi(\theta)]=\pi([\pi(\eta),\theta]+[\eta,\pi(\theta)])+
\pi\phi(\pi(\eta),\pi(\theta)),
\end{equation}
where $\phi$ is a closed 3-form on $V$.
It is known that if $\pi$ is a twisted Poisson structure
then the bracket $\{\eta,\theta\}:=[\pi(\eta),\theta]+[\eta,\pi(\theta)]+\phi(\pi(\eta),\pi(\theta))$ is a Lie bracket on the space of 1-forms.
In particular if $\phi=0$ then (\ref{tcybep}) coincides with (\ref{cybe2}).
Hence (non-twisted) Poisson structures are special examples
of twisted Poisson structures.
In this paper we construct a kind of quantum analogy
of (twisted-)Poisson geometry.
We here describe the main idea of this article.
\medskip\\
\indent
Carinena and coauthors \cite{CGM} introduced
a new operator identity which is
an associative version of classical Nijenhuis identity,
motivated by the study of Wigner problem.
The main result of \cite{CGM} is as follows.
Let $(\A,*)$ be an associative algebra equipped with
an operator $N:\A\to\A$.
They showed that
if $N$ satisfies (AN) below then
$x*_{1}y:=N(x)*y+x*N(y)-N(x*y)$ provides an associative
multiplication on $\A$ and the pair $(*,*_{1})$ becomes
a quantum bi-Hamiltonian system
in the sense of \cite{CGM}.
$$
N(x)*N(y)=N(N(x)*y+x*N(y))-NN(x*y).
\eqno{(AN)}
$$
The classical Nijenhuis condition is the Lie version of (AN),
$$
[N(x),N(y)]=N([N(x),y]+[x,N(y)])-NN([x,y]),
\eqno{(CN)}
$$
where $*$ in (AN) is replaced with a Lie bracket.
Hence the operator in (AN) is called an associative Nijenhuis operator.\\
\indent
We apply their idea to other operator identities, i.e.,
given an operator identity in the category of Lie algebras,
by the replacement of Lie brackets with associative multiplications,
we define an operator identity in the category of associative algebras.
We call such a formal functor a ``quantum analogy".
By using quantum analogy, we get a new operator identity,
that is an ``associative version" of the classical one.
We remark that quantum analogy is not unique,
for instance, the derivation rule in the category of associative
algebra is a quantum analogue
of the rule in the one of Lie algebras:
\begin{eqnarray*}
d(x*y)&=&dx*y+x*dy, \\
d[x,y]&=&[dx,y]+[x,dy].
\end{eqnarray*}
On the other hand, $d(x*y)=dx*y-dy*x$ is also
a quantum analogue, because $-[dy,x]=[x,dy]$.
However the latter is not interesting as analogy.
The associative Nijenhuis condition is considered
as a meaningful example of quantum analogues.
In general, it is difficult to construct meaningful examples.
We can find another interesting example
in Aguiar's works \cite{A0,A1,A2}.
Motivated by the study of infinitimal bi-algebras,
he introduced an associative Yang-Baxter equation
(shortly, AYBE):
$$
r_{13}*r_{12}-r_{12}*r_{23}+r_{23}*r_{13}=0,
\eqno{(AYBE)}
$$
which is an associative version of
classical Yang-Baxter equation (CYBE).
Here $r_{ij}$ are 3-power tensors and $*$ is an associative
multiplication.
The condition AYBE means that
if $r$ is a solution of AYBE then the Hochschild coboundary
$\partial r:\A\to\A\ot\A$ is a coassociative comultiplication.
This result is an analogy of a basic property of Lie bialgebras.
Hence (AYBE) can be seen as an integrability condition.
However we do not know yet an interesting geometrical application.\\
\indent
We introduce a quantum analogue of
the Poisson condition (\ref{cybe2}).
The Poisson bundle map $\pi:T^{*}V\to TV$ is defined
between the cotangent bundle and the tangent bundle.
Thus a quantum analogue of $\pi:T^{*}V\to TV$
is inappropriate for an endomorphism on $\A$.
The condition (\ref{cybe2}) is defined by using
a (non-symmetric) representation $TV\overset{rep}{\to}T^{*}V$.
So, in order to construct an associative version of (\ref{cybe2}),
we use the category of associative representations.
Let $M$ be an $\A$-bimodule, and let $\pi:M\to\A$ a linear map.
A quantum analogue of the Poisson condition
is defined to be the identity,
\begin{equation}\label{AP}
\pi(m)*\pi(n)=\pi(\pi(m)\cdot n+m\cdot\pi(n)),
\end{equation}
where $\cdot$ means the bimodule action.
We give several nontrivial examples of such operators
in the next section.
Especially, when $M=\A$ and $\cdot=*$, (\ref{AP}) is
reduced to a Rota-Baxter condition:
$$
\b(x)*\b(y)=\b(\b(x)*y+x*\b(y)),
\eqno{(RB)}
$$
where we changed $\pi$ for $\b$.
The $\b:\A\to\A$ is called a Rota-Baxter operator of ``weight zero",
or simply, Rota-Baxter operator.
In general, Rota-Baxter operators are defined by
the weighted formula, $\b(x)*\b(y)=\b(\b(x)*y+x*\b(y))-q\b(x*y)$,
where $q$ is a scalar (weight). The Rota-Baxter operators of $q=0$ are
special examples.
When an operator $\pi:M\to\A$ satisfies (\ref{AP}),
we call $\pi$ a \textbf{generalized Rota-Baxter operator},
or shortly, \textbf{GRB}-operator.
Remark that this generalization is only applied to
weight zero Rota-Baxter operators.
So, in the following, we assume the weight of Rota-Baxter
operator is zero.
In order to distinguish the usual Rota-Baxter operators
from the generalized ones,
we denote the usual Rota-Baxter operator by $\b$.
Mainly, the usual Rota-Baxter operators
have been studied in combinatorial theory
(\cite{Bax}, \cite{Rot1,Rot2}, \cite{Car}, \cite{FGK1,FGK2,FGK3}).
We do not discuss combinatorial problems,
because they are beyond our scope.\\
\indent
On the other hand, it is known that
Rota-Baxter operators $\b$ are closely related with
Loday's dendriform algebras in \cite{L1}.
Let $E$ be a module equipped with two binary
multiplications $\su$ and $\pr$.
$E$ is called a dendriform algebra,
if the following three axioms are satisfied.
\begin{eqnarray*}
(x\pr y)\pr z&=&x\pr(y\su z+y\pr z),\\
(x\su y) \pr z&=&x\su(y \pr z),\\
x\su (y\su z)&=&(x\su y+x\pr y)\su z,
\end{eqnarray*}
where $x,y,z\in E$.
He showed that the operad of dendriform algebras
is the Koszul dual of the one of associative dialgebras
and that a dendriform algebra characterizes an associative
multiplication, i.e., the sum of two multiplications,
$xy:=x\su y+x\pr y$, is associative.
In \cite{A1}, it was explained how one may associate a
dendriform algebra to any associative algebra equipped with a
Rota-Baxter operator, i.e., it was shown that if
$\b:\A\to\A$ is a Rota-Baxter operator then
the pair of multiplications $x\su y:=\b(x)y$ and $x\pr y:=x\b(y)$
is a dendriform structure on $\A$.
\medskip\\
\indent
In this letter we aim to make
a quantum analogue of Poisson geometry
by the generalized Rota-Baxter operators
and the theory of dendriform algebras.
Since (\ref{cybe2}) is equivalent with
the square zero equation (\ref{cybe})
and GRB-condition is an analogy of (\ref{cybe2}),
we hope to get a square zero equation
with respect to GRB-condition.
The bracket of (\ref{cybe}) is a Gerstenhaber-type bracket
defined by the derived bracket of
Kosmann-Schwarzbach \cite{Kos2}.
So, by using the derived bracket
associated with a canonical Gerstenhaber bracket,
we obtain the following square zero equation which
is equivalent with GRB-identity (Proposition \ref{bb=0} below).
$$
\frac{1}{2}[\hat{\pi},\hat{\pi}]_{\hat{\mu}}=0,
$$
where $\hat{\pi}$ is a natural extension of $\pi$
(Proposition \ref{defhat}).
Through Section 2, we will study
fundamental properties of generalized Rota-Baxter operators.
The generalized Rota-Baxter operators
inherit the basic properties from classical ones $\b:\A\to\A$.
For instance, an arbitrary GRB-operator
$\pi:M\to\A$ also induces a dendriform structure on $M$
by the natural manner.
However we will see several differences between
generalized Rota-Baxter operators and classical ones.
For instance, we can say that an arbitrary
dendriform structure is induced by a GRB-operator.
In other word, the notions
of ``Rota-Baxter" and ``dendriform" are unified by the generalization.
In Section 3, we define a quantum analogue of
the twisted Poisson condition (\ref{tcybep}).
We call a linear map $\pi:M\to\A$
a \textbf{twisted Rota-Baxter operator}, if
it is satisfying the condition,
$$
\pi(m)*\pi(n)=\pi(\pi(m)\cdot n+m\cdot\pi(n))+\pi\phi(\pi(m),\pi(n)),
\eqno{(TRB)}
$$
where $\phi$ is a Hochschild 2-cocycle in $C^{2}(\A,M)$.
A generalized Rota-Baxter operator is a special
twisted Rota-Baxter operator in which $\phi=0$.
Given an arbitrary 1-cochain $f$, one can define
the Hamiltonian vector field $X_{f}$ by usgin
the canonical Gerstenhaber bracket.
We will see that the Hamiltonian flow $exp(X_{f})$
is also well-defined in our construction (see Section 3, (\ref{hamflow})).
For an arbitrary 2-cochain $\Theta$,
its canonical transformation is defined
to be $exp(X_{f})(\Theta)$ by analogy with Poisson geometry.
By a standard argument of extension,
$$
0\longrightarrow M\longrightarrow\A\oplus M
\longrightarrow\A\longrightarrow 0,
$$
$\A\oplus M$ has an associative structure $\hat{\mu}+\hat{\phi}$,
where $\hat{\mu}$ is an extension of an associative structure of $\A$
and $\hat{\phi}$ is an extension of a 2-cocycle $\phi\in C^{2}(\A,M)$.
For a given 1-cochain $\pi$,
we will compute an explicit formula of
$exp(X_{\pi})(\hat{\mu}+\hat{\phi})$,
and from the result we obtain a geometrical characterization
of twisted Rota-Baxter operators
(Proposition \ref{bbepprop} and Corollary \ref{addexp}).
\medskip\\
\noindent
\textbf{Acknowledgements.}
The author wishes to thank very much the referees.
He is greatly indebted to them
for their numerous suggestions.
Especially,
a category theoretic approach in subsection 2.2
was given by referees.
Finally, he would like to thank very
much Professors Giuseppe Dito and Akira Yoshioka
for helpful comments and encouragement.

%%%%%%%%%%%%%%%%%%%%%%%%%%%%%%%%%%%%%%%%%%%%%%
\section{Rota-Baxter operators and dendriform algebras.}
%%%%%%%%%%%%%%%%%%%%%%%%%%%%%%%%%%%%%%%%%%%%%%

\noindent
\textbf{Notations and Standing assumptions.}
Throughout this article, $\A$ is an associative algebra over
a commutative ring $k$ and $M$ is an $\A$-bimodule.
The algebras are not necessarily unital.
Elements of $M$ and $\A$ will be denoted by
$m,n,l...$ and $a,b,c,...$ respectively.
And all actions will be denoted by dot ``$\cdot$".
In the following,
we omit the multiplication $*$ on $\A$,
$ab:=a*b$.

%%%%%%%%%%%%%%%%%%%%%%%%%%%%%%%%%%%%%%%%%%%%%%%%%%%
\subsection{Quantum analogue of Poisson condition.}
%%%%%%%%%%%%%%%%%%%%%%%%%%%%%%%%%%%%%%%%%%%%%%%%%%%

Let $\A$ be an associative algebra, and let $M$ an $\A$-bimodule,
and let $\pi:M\to\A$ a linear map.
We consider an operator identity,
\begin{equation}\label{defb}
\pi(m)*\pi(n)=\pi(\pi(m)\cdot n+m\cdot \pi(n)),
\end{equation}
where $m,n\in M$.
\begin{definition}\label{defbaxter}
We call $\pi:M\to \A$ a generalized Rota-Baxter operator,
or simply, GRB-operator, if it satisfies (\ref{defb}).
\end{definition}
A Rota-Baxter operator is a special example of GRB-operators
in which $M=\A$ and the bimodule structure is canonical.\\
\indent
Assume that $M^{\p}$ is an $\A$-bimodule.
One can easily check that
if $f:M^{\p}\to M$ is an $A$-bimodule morphism and $\pi:M\to\A$
is GRB then the composition $\pi\circ f:M^{\p}\to\A$
is also GRB. A classical Rota-Baxter
operator is known as an abstraction of the integral operator:
\begin{equation}\label{basicexample}
\b(f)(x):=\int^{x}_{0}dtf(t), \ \ f\in C^{0}([0,1]).
\end{equation}
\begin{example}\label{symanalog}
Given an invertible Hochschild 1-cocycle,
or invertible derivation $d:\A\to M$,
the inverse $\int:=d^{-1}$ is a GRB-operator.
For any $m,n\in M$, we have
$$
\int (m)\int (n)=\int d(\int (m)\int (n))=\int(m\cdot\int(n)+\int(m)\cdot n).
$$
\end{example}

\begin{example}\label{weyl}
Let $W\langle x,y\rangle$ be the Weyl algebra of two generators.
Here the commutation relation is $[x,y]:=1$.
For the normal basis $x^{i}y^{j}$, we define
$$
\int x^{i}y^{j}dy:=\frac{1}{j+1}x^{i}y^{j+1}.
$$
This integral operator is a classical Rota-Baxter operator.
\end{example}
\begin{proof}
When $j\neq 0$, we have
$$
\int [x,x^{i}y^{j}]dy=\int x^{i}[x,y^{j}]dy=
\int jx^{i}y^{j-1}dy=x^{i}y^{j}.
$$
For any $a=x^{i}y^{j}$, $b=x^{k}y^{l}$,
we have
$$
\int [x,\int ady\int bdy]dy=\int(\int adyb+a\int bdy)dy,
$$
where an identity $[x,\int(\cdot)dy]=id$ is used.
Since $\int ady\int bdy$ is a polynomial of the form
$C_{ij}x^{i}y^{j}$ and $j\neq 0$,
we have $\int [x,\int ady \int bdy]dy=\int ady\int bdy$.
\end{proof}
%\begin{remark}
%Similarly, $\int dx$ is also a classical Rota-Baxter operator,
%and $\int dx$ and $\int dy$ are commuting each other.
%(See \cite{A3} for the study of commuting pairs
%of Rota-Baxter operators.)
%\end{remark}
\begin{example}\label{newex1}
As the dual of the Rota-Baxter operator (\ref{basicexample}),
the derivation operator,
$$
\omega(x)\frac{d}{dx}:C^{1}([0,1])\to C^{0}([0,1]),
$$
becomes a generalized Rota-Baxter operator,
where $\omega(x)$ is a continuous function.
This example is followed from Proposition \ref{addob} below.
\end{example}

\begin{example}\label{weyl2}
As the dual of the integral operator in Example \ref{weyl},
the derivation operator,
$$
ad_{x}:W\langle x,y\rangle \to W\langle x,y\rangle,
\ \ a\mapsto [x,a]
$$
becomes GRB. This example is also followed
from Proposition \ref{addob} below.
\end{example}

We recall (AYBE) in Introduction:
\begin{equation}\label{rmatab}
\sum_{ij} a_{i}a_{j}\ot b^{j}\ot b^{i}=
\sum_{ij} a_{i}\ot b^{i}a_{j}\ot b^{j}-
\sum_{ij}a_{j}\ot a_{i}\ot b^{i}b^{j},
\end{equation}
where $r=\sum a_{i}\ot b_{i}\in\A\ot\A$.
We show that a skew symmetric solution of AYBE
is a GRB-operator. In the following, we omit $\sum_{ij}$.
\begin{example}
Let $r:=\sum_{i} a_{i}\ot b_{i}=a_{i}\ot b^{i}$
be a solution of AYBE.
We put $M:=\A^{*}:=\Hom_{k}(\A,k)$.
The $\A$-bimodule structure on $\A^{*}$ is
the usual one:
$(a\cdot f)(b):=f(ba)$ and $(f\cdot a)(b):=f(ab)$
for any $a,b\in\A$ and $f\in\A^{*}$.
We define a linear map $\ti{r}:\A^{*}\to \A$ by
$$
\ti{r}(f):=a_{i}f(b^{i}), \ \ f(b^{i})\in k.
$$
If $r$ is skew-symmetric then the $\ti{r}$ is GRB.
\end{example}
\begin{proof}
For any $f,g\in\A^{*}$, we have
$\ti{r}(f)\ti{r}(g)=a_{i}a_{j}g(b^{j})f(b^{i})$
and
\begin{multline*}
\ti{r}(\ti{r}(f)\cdot g)=a_{i}(\ti{r}(f)\cdot g)(b^{i})
=a_{i}g(b^{i}\ti{r}(f))=\\
=a_{i}g(b^{i}a_{j}f(b^{j}))=a_{i}g(b^{i}a_{j})f(b^{j}).
\end{multline*}
Assume that $r$ is skew symmetric.
Then we have $a_{i}f(b^{i})=-b_{i}f(a^{i})$ for any $f\in\A^{*}$.
By using the assumption, we have
\begin{multline*}
\ti{r}(f\cdot \ti{r}(g))=-\ti{r}(f\cdot g(a_{i})b^{i})=
-g(a_{i})\ti{r}(f\cdot b^{i})=\\
=-g(a_{i})a_{j}(f\cdot b^{i})(b^{j})=-a_{j}g(a_{i})f(b^{i}b^{j}).
\end{multline*}
We define a multilinear operation by
$(a\ot b\ot c)(g\ot f):=ag(b)f(c)$,
for any $g\ot f\in\A^{*}\ot\A^{*}$ and $a\ot b\ot c\in\A\ot\A\ot\A$.
Applying $g\ot f$ to (\ref{rmatab}),
the left hand-side is
$$
(a_{i}a_{j}\ot b^{j}\ot b^{i})(g\ot f)=
a_{i}a_{j}g(b^{j})f(b^{i})=\ti{r}(f)\ti{r}(g)
$$
and the right hand-side is
\begin{multline*}
(a_{i}\ot b^{i}a_{j}\ot b^{j})(g\ot f)
-(a_{j}\ot a_{i}\ot b^{i}b^{j})(g\ot f)=\\
=a_{i}g(b^{i}a_{j})f(b^{j})-a_{j}g(a_{i})f(b^{i}b^{j})=
\ti{r}(\ti{r}(f)\cdot g)+\ti{r}(f\cdot \ti{r}(g)).
\end{multline*}
From AYBE, we obtain the generalized
Rota-Baxter condition of $\ti{r}$.
\end{proof}
We study a useful representation of GRB-operators.
Set an algebra $\A\oplus_{0}M$ equipped with an
associative multiplication:
$$
(a,m)*(b,n):=(ab,a\cdot n+m\cdot b),
$$
where $(a,m),(b,n)\in\A\oplus M$.
Namely, $\A\oplus_{0}M$ is the algebra of semidirect product.
Let $\pi:M\to \A$ be a linear map. We denote
the graph of $\pi$ by $L_{\pi}$,
$$
L_{\pi}:=\{(\pi(m),m) \ | \ m\in M\}.
$$
In Poisson geometry,
it is well-known that $\pi$ is a Poisson structure on $V$
if and only if the graph of the bundle map $\pi:T^{*}V\to TV$
is a Lie algebroid which is called a Dirac structure (\cite{Cou}).
Namely, the Poisson condition is characterized by
certain algebraic properties of the graph.
As an analogy, we have Lemma \ref{thelemma} below.
We refer also to Freeman's early work \cite{Fre}.
He showed the lemma below in the classical cases.
\begin{lemma}\label{thelemma}
$\pi$ is a generalized Rota-Baxter operator
if and only if $L_{\pi}$ is a subalgebra of $\A\oplus_{0}M$.
\end{lemma}
\begin{proof}
For any $(\pi(m),m),(\pi(n),n)\in L_{\pi}$ we have
$$
(\pi(m),m)*(\pi(n),n)=(\pi(m)\pi(n),\pi(m)\cdot n+m\cdot\pi(n)).
$$
$\pi$ is a generalized Rota-Baxter if and only if
$(\pi(m)\pi(n),\pi(m)\cdot n+m\cdot\pi(n))$ is in $L_{\pi}$.
\end{proof}
The above lemma says that a generalized Rota-Baxter operator
is a good substructure of the trivial extension, $\A\oplus_{0}M$.
It is natural to ask what is a good substructure of an arbitrary
extension. We will give a solution in Section 3.\\
\indent
$L_{\pi}$ and $M$ are isomorphic as modules
by the identification $(\pi(m),m)\cong m$.
Thus if $\pi$ is a GRB-operator, i.e., $L_{\pi}$ is an associative
subalgebra of $\A\oplus_{0}M$ then $M$ is also an associative algebra.
\medskip\\
\indent
Given an arbitrary linear map $\pi:M\to\A$,
we define a lift of $\pi$, $\hat{\pi}$, as an endomorphism on $\A\oplus M$
by $\hat{\pi}(a,m):=(\pi(m),0)$.
\begin{proposition}\label{defhat}
$\pi:M\to\A$ is a generalized Rota-Baxter operator if and only if
the lift $\hat{\pi}(=:\b)$ is a classical Rota-Baxter operator
on $\A\oplus_{0}M$.
\end{proposition}
\begin{proof}
For any $(a,m),(b,n)\in\A\oplus_{0}M$, we have
$$
\hat{\pi}(a,m)*\hat{\pi}(b,n)=(\pi(m),0)*(\pi(n),0)=(\pi(m)\pi(n),0).
$$
On the other hand, we have
\begin{eqnarray*}
\hat{\pi}(\hat{\pi}(a,m)*(b,n)+(a,m)*\hat{\pi}(b,n))
&=&\hat{\pi}((\pi(m),0)*(b,n)+(a,m)*(\pi(n),0))\\
&=&\hat{\pi}((\pi(m)b,\pi(m)\cdot n)+(a\pi(n),m\cdot\pi(n)))\\
&=&(\pi(\pi(m)\cdot n+m\cdot\pi(n)),0).
\end{eqnarray*}
The proof of the proposition is completed.
\end{proof}
This proposition will be used in subsection 2.3.
\medskip\\
\indent
We recall dendriform algebras (\cite{L1}).
Let $E$ be a module equipped with two binary operations $\pr$ and $\su$.
$E$ is called a dendriform algebra, if the three axioms
(\ref{dend1}), (\ref{dend2}) and (\ref{dend3}) are satisfied.
For any $x,y,z\in E$,
\begin{eqnarray}
\label{dend1}(x\pr y)\pr z&=&x\pr(y\su z+y\pr z),\\
\label{dend2}(x\su y) \pr z&=&x\su(y \pr z),\\
\label{dend3}x\su (y\su z)&=&(x\su y+x\pr y)\su z,
\end{eqnarray}
In Lemma \ref{thelemma}
we saw that a GRB-operator $\pi:M\to\A$
induces an associative multiplication on $M$
via the isomorphism $M\cong L_{\pi}$.
This induced associative structure
comes from the dendriform algebra associated with $\pi$.
\begin{proposition}\label{phil}
(\cite{L1})
Given a dendriform algebra,
the sum of two multiplications $xy:=x\su y+x\pr y$
is an associative multiplication.
\end{proposition}
In \cite{A1} it was shown that if $\b:\A\to \A$ is
a classical Rota-Baxter operator
then $a\su b:=\b(a)b$ and $a\pr b:=a\b(b)$
define a dendriform algebra structure on $\A$.
For our generalized Rota-Baxter operators the same proposition holds.
\begin{proposition}\label{mainprop}
Let $\pi:M\to \A$ be a generalized Rota-Baxter operator.
Then $M$ becomes a dendriform algebra
by the multiplications $m\su n:=\pi(m)\cdot n$ and $m\pr n:=m\cdot \pi(n)$.
\end{proposition}
\begin{proof}
If $\pi$ is GRB then we have
\begin{eqnarray*}
m\su(n\su l)&=&(\pi(m)\pi(n))\cdot l\\
&=&\pi(\pi(m)\cdot n+m\cdot \pi(n))\cdot l\\
&=&(m\su n+m\pr n)\su l,
\end{eqnarray*}
where $l\in M$.
All other axioms are proved similarly.
\end{proof}
As an application of Proposition \ref{mainprop},
we will give Proposition \ref{addob} below.
When $\pi:M\to \A$ is GRB,
$M$ has the associative multiplication via
Propositions \ref{phil}, \ref{mainprop}.
We denote the associative algebra by $M_{ass}$.
\begin{lemma}\label{newaddlemma}
Under the assumptions above,
the module $\A$ becomes a $M_{ass}$-bimodule by the actions,
\begin{eqnarray*}
a\cdot_{\pi}m&:=&a\pi(m)-\pi(a\cdot m),\\
m\cdot_{\pi}a&:=&\pi(m)a-\pi(m\cdot a),
\end{eqnarray*}
where $\cdot_{\pi}$ means the $M_{ass}$-action on $\A$.
\end{lemma}
\begin{proof}
In Corollary \ref{mainprop2} below, we will show that the multiplication
$$
(a,m)*_{\pi}(b,n):=(a\cdot_{\pi}n+m\cdot_{\pi}b,mn)
$$
is associative on $\A\oplus M$, where $mn:=m\su n+m\pr n$.
The associativity gives the bimodule condition
of the action $\cdot_{\pi}$.
\end{proof}
The proposition below implies a duality
between integral operators and derivation operators.
\begin{proposition}\label{addob}
Let $\pi:M\to\A$ be a generalized Rota-Baxter operator,
and let $\Omega:\A\to M$ a derivation operator.
If $\Omega\pi(M)=z\cdot M$,
then $\Omega:\A\to M_{ass}$ is a generalized Rota-Baxter operator,
where $z$ is a central element in $Z(\A)$.
\end{proposition}
\begin{proof}
We have
\begin{eqnarray*}
\Omega(a)\cdot_{\pi}b+a\cdot_{\pi}\Omega(b)
&=&\pi(\Omega(a))b-\pi(\Omega(a)\cdot b)+a\pi(\Omega(b))
-\pi(a\cdot \Omega(b))\\
&=&\pi(\Omega(a))b+a\pi(\Omega(b))-\pi(\Omega(a)\cdot b+a\cdot \Omega(b))\\
&=&\pi(\Omega(a))b+a\pi(\Omega(b))-\pi\Omega(ab).
\end{eqnarray*}
Here the derivation condition of $\Omega$ is used.
Applying $\Omega$ on the both side,
we have
\begin{eqnarray*}
\Omega(\Omega(a)\cdot_{\pi}b+a\cdot_{\pi}\Omega(b))&=&
\Omega(\pi(\Omega(a))b+a\pi(\Omega(b))-\pi\Omega(ab))\\
&=&\Omega(\pi(\Omega(a))b+a\pi(\Omega(b)))-z\cdot\Omega(ab)\\
&=&z\Omega(a)\cdot b+\pi(\Omega(a))\cdot\Omega(b)+
\Omega(a)\cdot\pi(\Omega(b))+za\cdot\Omega(b)-z\cdot\Omega(ab)\\
&=&\pi(\Omega(a))\cdot\Omega(b)+\Omega(a)\cdot\pi(\Omega(b))\\
&=&\Omega(a)\su \Omega(b)+\Omega(a)\pr \Omega(b)=\Omega(a)\Omega(b).
\end{eqnarray*}
Here $\Omega\pi=z$ is used.
\end{proof}
We remember Examples \ref{newex1} and \ref{weyl2}.
In the both examples,
the assumption of the proposition above
holds. Hence the derivation operators of the two examples
become generalized Rota-Baxter operators.
\begin{remark}
Under the assumptions of Proposition \ref{addob},
one can show that $\pi\Omega:\A\to\A$ is an associative
Nijenhuis operator satisfying (AN) in Introduction.
This implies a relationship between generalized Rota-Baxter operators
and associative Nijenhuis operators.
We will study the relation in \cite{U}.
\end{remark}

%%%%%%%%%%%%%%%%%%%%%%%%%%%%%%%%%%%%%%%%%%%%%%%%%%%%%%%%%%%
\subsection{Rota-Baxter operators via Dendriform algebras.}
%%%%%%%%%%%%%%%%%%%%%%%%%%%%%%%%%%%%%%%%%%%%%%%%%%%%%%%%%%%

We consider the converse of Proposition \ref{mainprop}.
Given a dendriform algebra $E$,
$E$ is an associative algebra by Proposition \ref{phil}.
We denote the associated associative algebra by $E_{ass}$.
One can easily check that $E$ is an $E_{ass}$-bimodule
by $e\cdot x:=e\su x$ and $x\cdot e:=x\pr e$,
where $e\in E_{ass}$ and $x\in E$.
Under this setting, the identity map $1:E\to E_{ass}$
is a {\em generalized} Rota-Baxter operator
and the corresponding dendriform algebra is
the original one.
Hence all dendriform algebras are induced by generalized
Rota-Baxter operators.
However this correspondence between
generalized Rota-Baxter operators and dendriform algebras
is not bijective. We will discuss the correspondence as follow.
\medskip\\
\indent
Let $\Dend$ be the category of dendriform algebras.
The objects of $\Dend$ are dendriform algebras and
the morphisms are the obvious ones.
In addition we consider the category of GRB-operators,
it is denoted by $\Bax$.
The objects of $\Bax$ are GRB-operators $\pi:M\to\A$
and the morphisms are commutative diagrams of linear maps:
\begin{equation}\label{defmor}
\begin{CD}
M@>\psi_{1}>>M^{\p} \\
@V\pi VV @V\gamma VV \\
\A@>\psi_{0}>>\A^{\p}
\end{CD}
\end{equation}
such that $\psi_{1}(a\cdot m)=\psi_{0}(a)\cdot\psi_{1}(m)$
and $\psi_{1}(m\cdot a)=\psi_{1}(m)\cdot\psi_{0}(a)$,
where $\gamma$ is a GRB-operator
and $\psi_{0}$ is not necessarily an algebra homomorphism.
Proposition \ref{mainprop}
defines a functor $G:\Bax\to\Dend$.
Conversely, for a given dendriform algebra $E$,
the GRB-operator of identity map $1:E\to E_{ass}$ is
the result of a functor $F:\Dend\to\Bax$.
The composition $GF$ is clearly identity,
and one can check an adjoint relation of the pair $(F,G)$:
\begin{equation}\label{adre}
\Hom_{\Bax}(F(E),\pi)\cong\Hom_{\Dend}(E,G(\pi)),
\end{equation}
where $\pi$ means an object $\pi:M\to \A$.
This isomorphism sends an arbitrary morphism
$(\psi:E\to M)\in\Hom_{\Dend}(E,G(\pi))$
to the following morphism in $\Hom_{\Bax}(F(E),\pi)$.
\begin{equation}\label{adre2}
\begin{CD}
E@>\psi>>M \\
@V1VV @V\pi VV \\
E_{ass}@>\pi\circ\psi>>\A.
\end{CD}
\end{equation}
%%%%%%%%%%%%%%%%%%%%%%%%%%%%%%%%%
\subsection{Structure equation.}
%%%%%%%%%%%%%%%%%%%%%%%%%%%%%%%%%
In this subsection we give the GRB version
of the Poisson condition $[\pi,\pi]=0$.
In the following, we assume that the base ring $k$
has $1/2$.
\medskip\\
\indent
We denote the associative structure of $\A$
by $\mu:\A\ot\A\to\A$. We define $\hat{\mu}$
as the associative multiplication of $\A\oplus_{0}M$ by
$\hat{\mu}(a,b)=ab$,
$\hat{\mu}(a,n)=a\cdot n$,
$\hat{\mu}(m,b)=m\cdot b$ and
$\hat{\mu}(m,n)=0$.
We set the space:
$$
G(\A\oplus_{0}M):=
\bigoplus_{n\ge 1}\Hom((\A\oplus_{0}M)^{\ot n},\A\oplus_{0}M).
$$
The space $G(\A\oplus_{0}M)$ has a bracket product $[\cdot,\cdot]$
which is called a Gerstenhaber-bracket, or shortly G-bracket
(See Appendix for the definition of G-brackets).
It is well-known that $S\in\Hom((\A\oplus_{0}M)^{\ot 2},\A\oplus_{0}M)$
is a solution of $[S,S]=0$ if and only if
$S$ is an associative multiplication on $\A\oplus_{0}M$.
Hence we have $[\hat{\mu},\hat{\mu}]=0$.
By the graded Jacobi rule of G-bracket,
$d_{\hat{\mu}}:=[\hat{\mu},\cdot]$ becomes a square-zero
derivation of degree $+1$.
By using the derivation, we define the second bracket below.
For any  $f,g\in G(\A\oplus_{0}M)$,
$$
[f,g]_{\hat{\mu}}:=[d_{\hat{\mu}}f,g]=[[\hat{\mu},f],g].
$$
The bracket, $[\cdot,\cdot]_{\hat{\mu}}$, is called a derived bracket
(see \cite{Kos2}).
The derived bracket is not graded commutative,
but it satisfies a graded Leibniz rule.\\
\indent
The derivation of $\hat{\pi}$ by $\hat{\mu}$
has the form,
$$
[\hat{\mu},\hat{\pi}]=\hat{\mu}(\hat{\pi}\ot id)+\hat{\mu}(id\ot\hat{\pi})
-\hat{\pi}\c\hat{\mu}.
$$
where $\hat{\pi}$ is defined in Proposition \ref{defhat}.
The following lemma is the key.
\begin{lemma}
$\hat{\pi}\c\hat{\pi}=0$.
\end{lemma}
From this lemma we obtain
$$
\frac{1}{2}[\hat{\pi},\hat{\pi}]_{\hat{\mu}}=
\hat{\mu}(\hat{\pi}\ot\hat{\pi})-\hat{\pi}\c\hat{\mu}(\hat{\pi}\ot id)
-\hat{\pi}\c\hat{\mu}(id\ot\hat{\pi}),
$$
or explicitly, for any $(a,m),(b,n)\in\A\oplus_{0}M$,
\begin{multline}\label{mmb2}
\frac{1}{2}[\hat{\pi},\hat{\pi}]_{\hat{\mu}}((a,m),(b,n))=\\
\hat{\mu}((\pi(m),0),(\pi(n),0))-\hat{\pi}\c\hat{\mu}((\pi(m),0),(b,n))
-\hat{\pi}\c\hat{\mu}((a,m),(\pi(n),0))=\\
\pi(m)\pi(n)-\pi(\pi(m)\cdot n)-\pi(m\cdot\pi(n)).
\end{multline}
From (\ref{mmb2}), we obtain
the structure equation for GRB-operator.
\begin{proposition}\label{bb=0}
A linear map $\pi:M\to\A$ is a GRB-operator
if and only if it is a solution of
\begin{equation*}
\frac{1}{2}[\hat{\pi},\hat{\pi}]_{\hat{\mu}}=0.
\end{equation*}
\end{proposition}
From the graded Jacobi rule of G-bracket,
we obtain corollaries below.
\begin{corollary}\label{mainprop2}
If $\pi$ is GRB then
$[\hat{\mu},\hat{\pi}]$ is an associative structure on $\A\oplus M$,
and the associative multiplication has the form,
$$
(a,m)*_{\pi}(b,n):=
[\hat{\mu},\hat{\pi}]((a,m),(b,n))=(a\cdot_{\pi}n+m\cdot_{\pi}b,mn)
$$
where $mn=\pi(m)\cdot n+m\cdot\pi(n)$ and
$\cdot_{\pi}$ was defined in Lemma \ref{newaddlemma}.
\end{corollary}
\begin{proof}
Since $[[\hat{\mu},\hat{\pi}],[\hat{\mu},\hat{\pi}]]=0$,
it is an associative structure.
For any $a,b\in\A$ and $m,n\in M$, we obtain
$$
[\hat{\mu},\hat{\pi}](a,m)=
\hat{\mu}(\hat{\pi}(a)\ot m)+\hat{\mu}(a\ot\hat{\pi}(m))
-\hat{\pi}\c\hat{\mu}(a,m)=
a\pi(m)-\pi(a\cdot m)
$$
and $[\hat{\mu},\hat{\pi}](n,b)=\pi(n)b-\pi(n\cdot b)$.
In the same way, we have
$[\hat{\mu},\hat{\pi}](m,n)=\pi(m)\cdot n+m\cdot\pi(n)$
and $[\hat{\mu},\hat{\pi}](a,b)=0$.
\end{proof}
By above corollary, the proof of Lemma \ref{newaddlemma} is completed.
\begin{corollary}
$[\hat{\pi},\cdot]_{\hat{\mu}}$
is a square zero derivation of degree $+1$
for the derived bracket.
\end{corollary}
When $\pi$ is GRB, $\hat{\mu}+[\hat{\mu},\hat{\pi}]$
is an associative structure, i.e.,
the pair of associative structures,
$(\hat{\mu},[\hat{\mu},\hat{\pi}])$, is compatible.
This will be extended in Corollary \ref{addexp} below.

%%%%%%%%%%%%%%%%%%%%%%%%%%%%%%%%%%%%%%%%%%%%%%%%%%%%%%%%
\section{Twisted Rota-Baxter operators and NS-algebras.}
%%%%%%%%%%%%%%%%%%%%%%%%%%%%%%%%%%%%%%%%%%%%%%%%%%%%%%%%

In this section,
we will construct the twisted version of Section 2.
\medskip\\
\indent
The semi-direct product algebra appears as
the trivial extension of $\A$ by $M$:
$$
\begin{CD}
0@>>>M@>>>\A\oplus_{0}M@>>>\A@>>>0. \ \ \text{(exact)}
\end{CD}
$$
In general, given an abelian extension,
an associative multiplication on $\A\oplus M$
has the form (\cite{L4}),
\begin{equation}\label{twistedm}
(a,m)*(b,n):=(ab,a\cdot n+m\cdot b+\phi(a,b)),
\end{equation}
where $\phi$ is a Hochschild 2-cocycle in $C^{2}(\A,M)$, i.e.,
$\phi$ is satisfying a cocycle condition of Hochschild:
$$
0=\partial\phi(a,b,c)=a\cdot\phi(b,c)-\phi(ab,c)+\phi(a,bc)-\phi(a,b)\cdot c,
$$
where $\partial$ is Hochschild's coboundary map and $a,b,c\in\A$.
We denote an associative algebra $\A\oplus M$ equipped with
the twisted multiplication (\ref{twistedm}) by $\A\oplus_{\phi}M$.
\begin{definition}\label{deftwb}
Let $\pi:M\to \A$ be a linear map.
We call $\pi$ a \textbf{twisted Rota-Baxter operator},
or simply \textbf{$\phi$-Rota-Baxter operator}, if
\begin{equation}\label{defbt}
\pi(m)\pi(n)=\pi(\pi(m)\cdot n+m\cdot \pi(n))+\pi\phi(\pi(m),\pi(n)).
\end{equation}
is satisfied for any $m,n\in M$.
\end{definition}
Similarly to Lemma \ref{thelemma}, we consider the graph of $\pi$.
It is shown that the graph $L_{\pi}$ is a subalgebra
of $\A\oplus_{\phi}M$ if and only if
$\pi$ is a $\phi$-Rota-Baxter operator.
From the isomorphism $L_{\pi}\cong M$,
a $\phi$-Rota-Baxter operator
induces an associative multiplication on $M$.
The induced multiplication on $M$ has the form,
\begin{equation}\label{addtimes}
m\times n:=\pi(m)\cdot n+m\cdot \pi(n)+\phi(\pi(m),\pi(n)).
\end{equation}
It is obvious that $\pi$ is an algebra homomorphism:
$\pi(m\times n)=\pi(m)\pi(n)$.
\begin{example}\label{symplectic2}
Let $\omega:M\to\A$ be an invertible 1-cochain.
Then the inverse $\omega^{-1}$ is a twisted
Rota-Baxter operator, and in this case $\phi:=-\partial \omega$.
\end{example}
\begin{proof}
We put $\pi:=\omega^{-1}$.
The condition (\ref{defbt}) is equal to
$$
\omega(\pi(m)\pi(n))=\pi(m)\cdot n+m\cdot \pi(n)+\phi(\pi(m),\pi(n)).
$$
This is the same as
$-\pi(m)\cdot n+\omega(\pi(m)\pi(n))-m\cdot \pi(n)=\phi(\pi(m),\pi(n))$.
Since the coboundary of $\omega$ is defined by
$\partial\omega(a,b)=a\cdot\omega(b)-\omega(ab)+\omega(a)\cdot b$,
we obtain $\partial\omega(\pi(m),\pi(n))=-\phi(\pi(m),\pi(n))$.
\end{proof}

\begin{example}
Assume that $\A$ is unital.
Let $f:M\to\A$ be a $\A$-linear surjection.
We fix an element $e\in M$ such that $f(e)=1_{\A}$,
where $1_{\A}$ is the unit element of $\A$.
A cochain $\phi(a,b):=-a\cdot e\cdot b$ is a 2-cocycle.
$f$ is twisted Rota-Baxter:
\begin{multline*}
f(f(m)\cdot n+m\cdot f(n))-f(f(m)\cdot e\cdot f(n))=\\
=f(m)f(n)+f(m)f(n)-f(m)1_{\A}f(n)=f(m)f(n).
\end{multline*}
For instance, $M=\A\ot\A$, $f=\mu$ and $e=x\ot x^{-1}$,
where $\mu$ is the associative multiplication
and $x$ is an invertible element.
\end{example}

\begin{example}\label{exampleaota}
$\A$ is not necessarily unital.
An associative multiplication $\mu:\A\ot\A\to\A$, $\mu(a\ot b)=ab$
is a twisted Rota-Baxter operator with a cocycle $\phi(a,b):=-a\ot b$.
\end{example}

\begin{example}
A linear endomorphism $R:\A\to\A$ is called a \textbf{Reynolds operator}
(\cite{Rot2}), if the condition
$$
R(a)R(b)=R(R(a)b+aR(b))-R(R(a)R(b))
$$
is satisfied for any $a,b\in\A$. The last term $-R(a)R(b)$
is the associative multiplication $-\mu(R(a),R(b))$
which is the Hochschild 2-cocycle.
Thus each Reynolds operator can be seen as
twisted Rota-Baxter operator.
\end{example}
In Proposition \ref{mainprop},
we saw that GRB-operators induce dendriform
algebra structures.
We show a similar result with respect to twisted Rota-Baxter operators.
First we recall NS-algebras of Leroux \cite{Ler}.
\begin{definition}\label{tdend}
Let $\T$ be a $k$-module equipped with three binary multiplications
$\su$, $\pr$ and $\vee$. $\T$ is called a \textbf{NS-algebra}, if
the following 4 axioms are satisfied. For any $x,y,z\in\T$:
\begin{eqnarray}
\label{td1}(x\pr y)\pr z&=&x\pr(y\su z+y\pr z+y\vee z),\\
\label{td2}(x\su y)\pr z&=&x\su(y\pr z),\\
\label{td3}x\su(y\su z)&=&(x\su y+x\pr y+x\vee y)\su z,
\end{eqnarray}
and
\begin{equation}\label{td4}
x\su(y\vee z)-(x\times y)\vee z+x\vee(y\times z)-(x\vee y)\pr z=0,
\end{equation}
where $\times$ is defined as $x\times y:=x\su y+x\pr y+x\vee y$.
\end{definition}
We recall the associative Nijenhuis condition (AN)
in Introduction. In \cite{Ler},
it was shown that if $N$ is an associative Nijenhuis
operator then the multiplications
$x\su y:=N(x)y$, $x\pr y:=xN(y)$ and $x\vee y:=-N(xy)$
satisfy the axioms of NS-algebras.
The basic property of usual dendriform algebras is
satisfied on NS-algebras.
\begin{proposition}\label{timesass}
(\cite{Ler})
Let $\T$ be a NS-algebra.
Then the multiplication $\times$ above
is associative. 
\end{proposition}
We show that the quantum algebra, or nonclassical algebra,
associated with a twisted Rota-Baxter operator
is a NS-algebra.
\begin{proposition}\label{maintheorem2}
If $\pi:M\to \A$ is a $\phi$-Rota-Baxter operator
then $M$ is a NS-algebra
by the following three multiplications.
$$
m\su n:=\pi(m)\cdot n, \ \ m\pr n:=m\cdot\pi(n),
\ \ m\vee n:=\phi(\pi(m),\pi(n))
$$
\end{proposition}
\begin{proof}
It is easy to check the conditions (\ref{td1})-(\ref{td3}).
We show (\ref{td4}). Since $\phi$ is a Hochschild cocycle,
we have the cocycle condition:
\begin{multline*}
\partial\phi(\pi(m),\pi(n),\pi(l))=
\pi(m)\cdot\phi(\pi(n),\pi(l))-\phi(\pi(m)\pi(n),\pi(l))+\\
\phi(\pi(m),\pi(n)\pi(l))-\phi(\pi(m),\pi(n))\cdot \pi(l)=0,
\end{multline*}
where $m,n,l\in M$.
From the definitions of $\pr$, $\su$ and $\vee$,
we have
\begin{equation}\label{prfbd}
m\su(n\vee l)-\phi(\pi(m)\pi(n),\pi(l))+
\phi(\pi(m),\pi(n)\pi(l))-(m\vee n)\pr l=0.
\end{equation}
On the other hand, we have
\begin{eqnarray*}
\pi(m)\pi(n)&=&\pi(\pi(m)\cdot n+m\cdot \pi(n)+\phi(\pi(m),\pi(n))\\
&=&\pi(m\su n+m\pr n+m\vee n)\\
&=&\pi(m\times n).
\end{eqnarray*}
Hence we obtain
$$
\phi(\pi(m)\pi(n),\pi(l))=\phi(\pi(m\times n),\pi(l))=(m\times n)\vee l.
$$
In the same way, we have
$\phi(\pi(m),\pi(n)\pi(l))=\phi(\pi(m),\pi(n\times l))=m\vee(n\times l)$.
This says that (\ref{prfbd}) is equal to (\ref{td4}).
\end{proof}

We study a relation between twisted Rota-Baxter operators
and NS-algebras.
Let $\T$ be a NS-algebra.
We denote the associative algebra $(\T,\times)$ by $\T_{ass}$.
A $\T_{ass}$-bimodule structure on $\T$ is well-defined
by $t\cdot x:=t\su x$ and $x\cdot t:=x\pr t$,
where $t\in\T_{ass}$.
\begin{lemma}\label{bco}
For any $x,y\in\T_{ass}$, $\Phi(x,y):=x\vee y$ is a Hochschild cocycle
in $C^{2}(\T_{ass},\T)$.
\end{lemma}
\begin{proof}
The cocycle condition of $\Phi$
is the same as (\ref{td4}).
\end{proof}
We denote the categories of twisted Rota-Baxter
operators and NS-algebra by
$\phiBax$ and $\NS$, respectively.
The objects of $\NS$ are NS-algebras and
the morphisms are the obvious ones.
The objects and morphisms of $\phiBax$
are defined by the following commutative diagram.
$$
\begin{CD}
\A\ot\A@>\psi_{0}\ot\psi_{0}>>\A^{\p}\ot\A^{\p} \\
@V{\phi}VV @V\phi^{\p}VV \\
M@>\psi_{1}>>M^{\p} \\
@V{\pi}VV @V{\gamma}VV \\
\A@>\psi_{0}>>\A^{\p}.
\end{CD}
$$
The compatibility conditions between $\psi_{0}$ and $\psi_{1}$
are defined by the same way as subsection 2.2.
From Proposition \ref{maintheorem2},
we have a functor $G$ from $\phiBax$ to $\NS$.
On the other hand,
by Lemma \ref{bco},
if $\T$ is a NS-algebra then
the identity map $1:\T\to\T_{ass}$ is a twisted Rota-Baxter
operator with the 2-cocycle $\Phi(x,y)=x\vee y$.
This defines a functor $F:\NS\to\phiBax$.
We obtain an adjoint relation:
$$
\Hom_{\phiBax}(F(\T),(\pi,\phi))\cong\Hom_{\NS}(\T,G(\pi,\phi)),
$$
where $(\pi,\phi)$ is an object in $\phiBax$.
\medskip\\
\indent
Recall Proposition \ref{bb=0}.
We consider a structure equation for twisted Rota-Baxter operators.
For any linear map $\pi:M\to\A$, we have
$$
\frac{1}{2}[\hat{\pi},\hat{\pi}]_{\hat{\mu}}=
\hat{\mu}(\hat{\pi}\ot\hat{\pi})-\hat{\pi}\c\hat{\mu}(\hat{\pi}\ot id)
-\hat{\pi}\c\hat{\mu}(id\ot\hat{\pi})
$$
where $\hat{\mu}$ was defined in subsection 2.3.
On the other hand, $\pi$ is a $\phi$-Rota-Baxter operator
if and only if $\hat{\pi}$ satisfies
\begin{eqnarray}\label{add226}
\hat{\mu}(\hat{\pi}\ot\hat{\pi})
-\hat{\pi}\c\hat{\mu}(\hat{\pi}\ot id)
-\hat{\pi}\c\hat{\mu}(id\ot\hat{\pi})
-\hat{\pi}\c\hat{\phi}(\hat{\pi}\ot\hat{\pi})=0,
\end{eqnarray}
where $\hat{\phi}(a,b):=\phi(a,b)$ and $\hat{\phi}=0$ all other cases.
We recall the notion of twisted Poisson structure (\ref{tcybep})
in Introduction. It is known that the twisted Poisson condition
is equivalent with a modified Maurer-Cartan equation:
\begin{eqnarray}\label{tpstruct}
\frac{1}{2}[\pi,\pi]=-\frac{1}{6}\{\{\{\phi,\pi\},\pi\},\pi\},
\end{eqnarray}
where $\{\cdot,\cdot\}$ is a certain graded Poisson bracket.
We give a twisted Rota-Baxter version of (\ref{tpstruct}).
\begin{proposition}\label{bbepprop}
We assume $\mathbb{Q}\subset k$.
A linear map $\pi:M\to\A$ is a $\phi$-Rota-Baxter operator
if and only if $\pi$ is a solution of
\begin{equation*}
\frac{1}{2}[\hat{\pi},\hat{\pi}]_{\hat{\mu}}=
-\frac{1}{6}[[[\hat{\phi},\hat{\pi}],\hat{\pi}],\hat{\pi}].
\end{equation*}
\end{proposition}
\begin{proof}
\begin{eqnarray*}
\frac{1}{2}[[[\hat{\phi},\hat{\pi}],\hat{\pi}],\hat{\pi}]
&=&[\hat{\phi}(\hat{\pi}\ot\hat{\pi})-\hat{\pi}\c\hat{\phi}(\hat{\pi}\ot 1)-
\hat{\pi}\c\hat{\phi}(1\ot\hat{\pi}),\hat{\pi}] \\
&=&[\hat{\phi}(\hat{\pi}\ot\hat{\pi}),\hat{\pi}]-
[\hat{\pi}\c\hat{\phi}(\hat{\pi}\ot 1),\hat{\pi}]-
[\hat{\pi}\c\hat{\phi}(1\ot\hat{\pi}),\hat{\pi}] \\
&=&-\hat{\pi}\c\hat{\phi}(\hat{\pi}\ot\hat{\pi})
-\hat{\pi}\c\hat{\phi}(\hat{\pi}\ot\hat{\pi})
-\hat{\pi}\c\hat{\phi}(\hat{\pi}\ot\hat{\pi}) \\
&=&-3\hat{\pi}\c\hat{\phi}(\hat{\pi}\ot\hat{\pi}).
\end{eqnarray*}
Thus we have $\hat{\pi}\c\hat{\phi}(\hat{\pi}\ot\hat{\pi})=
-\frac{1}{6}[[[\hat{\phi},\hat{\pi}],\hat{\pi}],\hat{\pi}]$.
From (\ref{add226}), we obtain the desired result.
\end{proof}

Let $\pi:M\to\A$ be a linear map, not necessarily (twisted-)Rota-Baxter.
We define a Hamiltonian vector field of $\pi$
by the derivation, $X_{\pi}(\cdot):=[\cdot,\hat{\pi}]$.
The associated Hamiltonian flow is, by definition,
\begin{equation}\label{hamflow}
exp(X_{\pi}):=1+X_{\pi}+\frac{1}{2!}X_{\pi}^{2}+
\frac{1}{3!}X_{\pi}^{3}+...,
\end{equation}
where the series is convergent, because $\hat{\pi}\c\hat{\pi}=0$.
\begin{remark}
One can directly check that
$exp(X_{\pi})(\hat{\mu}+\hat{\phi})$
is again an associative structure
and that it is isomorphic to $\hat{\mu}+\hat{\phi}$ on $\A\oplus M$.
\end{remark}
We should determine the structure of
$exp(X_{\pi})(\hat{\mu}+\hat{\phi})$.
We have
$$
\frac{1}{2!}X_{\pi}^{2}(\hat{\mu}+\hat{\phi})
=\frac{1}{2!}[[\hat{\mu}+\hat{\phi},\hat{\pi}],\hat{\pi}]
=\frac{1}{2!}[\hat{\pi},\hat{\pi}]_{\hat{\mu}}
+\frac{1}{2!}[[\hat{\phi},\hat{\pi}],\hat{\pi}]
$$
and
$$
\frac{1}{3!}X_{\pi}^{3}(\hat{\mu}+\hat{\phi})
=\frac{1}{3!}[[[\hat{\mu}+\hat{\phi},\hat{\pi}],\hat{\pi}],\hat{\pi}]
=\frac{1}{3!}[[[\hat{\phi},\hat{\pi}],\hat{\pi}],\hat{\pi}],
$$
where $[[[\hat{\mu},\hat{\pi}],\hat{\pi}],\hat{\pi}]=0$ is used.
Hence we obtain
$$
exp(X_{\pi})(\hat{\mu}+\hat{\phi})=
\hat{\mu}+\hat{\phi}+[\hat{\mu}+\hat{\phi},\hat{\pi}]
+\frac{1}{2}[\hat{\pi},\hat{\pi}]_{\hat{\mu}}
+\frac{1}{2}[[\hat{\phi},\hat{\pi}],\hat{\pi}]
+\frac{1}{6}[[[\hat{\phi},\hat{\pi}],\hat{\pi}],\hat{\pi}],
$$
where $X^{I}_{\pi}(\hat{\mu}+\hat{\phi})$ are all zero for any $4\le I$,
because $\hat{\pi}\c\hat{\pi}=0$.
As a corollary of Proposition \ref{bbepprop}, we give
\begin{corollary}\label{addexp}
$\pi$ is a $\phi$-Rota-Baxter operator if and only if
$$
exp(X_{\pi})(\hat{\mu}+\hat{\phi})=
\hat{\mu}+\hat{\phi}+[\hat{\mu}+\hat{\phi},\hat{\pi}]
+\frac{1}{2}[[\hat{\phi},\hat{\pi}],\hat{\pi}].
$$
Then $exp(X_{\pi})(\hat{\mu}+\hat{\phi})$ defines an associative
multiplication on $M$
which has the form (\ref{addtimes}).
\end{corollary}
\begin{proof}
For any $m,n\in M$, we have
$$
exp(X_{\pi})(\hat{\mu}+\hat{\phi})(m,n)
=[\hat{\mu},\hat{\pi}](m,n)
+\frac{1}{2}[[\hat{\phi},\hat{\pi}],\hat{\pi}](m,n),
$$
where all other terms are zero.
Simply, we obtain
$$
[\hat{\mu},\hat{\pi}](m,n)
+\frac{1}{2}[[\hat{\phi},\hat{\pi}],\hat{\pi}](m,n)=
\pi(m)\cdot n+m\cdot \pi(n)+\phi(\pi(m),\pi(n)).
$$
\end{proof}
%%%%%%%%%%%%%%%%%%%%
\section{Appendix.}
%%%%%%%%%%%%%%%%%%%%
Let $B$ be a $k$-module.
We set the space of multilinear maps,
$$
G(B):=\bigoplus_{n\ge 1}\Hom_{k}(B^{\ot n},B).
$$
The degree of $f\in G(B)$ is $m$, if $f$ is in $\Hom_{k}(B^{\ot m},B)$.
For any $f\in\Hom_{k}(B^{\ot m},B)$ and
$g\in\Hom_{k}(B^{\ot n},B)$, define a binary product $\bar{\c}$:
$$
f\bar{\c}g:=\sum^{m}_{i=1}(-1)^{(i-1)(n-1)}f\c_{i}g,
$$
where $\c_{i}$ is the composition of maps defined by
$$
f\c_{i}g(b_{1},...,b_{m+n})=
f(b_{1},...,b_{i-1},g(b_{i},...,b_{i+n-1}),b_{i+n}...,b_{m+n}).
$$
The degree of $f\bar{\c}g$ is $m+n-1$.
The G-bracket on $G(B)$ is by definition a graded commutator:
$$
[f,g]:=f\bar{\c}g-(-1)^{(m-1)(n-1)}g\bar{\c}f.
$$
We recall two fundamental identities:
$[f,g]=-(-1)^{(m-1)(n-1)}[g,f]$ and
\begin{multline*}
(-1)^{(m-1)(l-1)}[[f,g],h]+(-1)^{(l-1)(n-1)}[[h,f],g]+\\
(-1)^{(n-1)(m-1)}[[g,h],f]=0,
\end{multline*}
where the degree of $h$ is $l$.
The above graded Jacobi rule is equivalent with the following
graded Leibniz rule.
$$
[f,[g,h]]=[[f,g],h]+(-1)^{(m-1)(n-1)}[g,[f,h]].
$$

%%%%%%%%%%%%%%%%%%%%%%%%%%%%%%%%%%%%%%%%%%%%%%%%%%%%%

\ \\
\noindent
Kyousuke UCHINO, Post doctoral student,
Tokyo University of Science,
Shinjyuku Wakamiya 26, Tokyo Japan.
e-mail: K\underline{ }Uchino[at]oct.rikadai.jp\\


\begin{thebibliography}{}

\bibitem{A0}
M. Aguiar.
Infinitesimal Hopf algebras.
New trends in Hopf algebra theory.
267 (1999), 1--29.
\bibitem{A1}
M. Aguiar.
Pre-Poisson algebras.
Lett. Math. Phys.
54 (2000), no 4, 263--277.
\bibitem{A2}
M. Aguiar.
On the associative analog of Lie bialgebras.
J. Algebra.
244 (2001), no 2, 492--532.
%\bibitem{A3}
%M. Aguiar and J-L. Loday.
%Quadri-algebras.
%J. Pure Applied Algebra.
%191 (2004), 205--221.
\bibitem{Bax}
G. Baxter.
An analytic problem whose solution follows from a simple algebraic identity.
Pacific J. Math.
10 (1960), 731--742.
\bibitem{CGM}
J.F. Carinena, J. Grabowski and G. Marmo.
Quantum Bi-Hamiltonian Systems.
Int. J. Mod. Phys.
A15 (2000), 4797-4810.
\bibitem{Car}
P. Cartier.
On the structure of free Baxter algebras.
Advances in Math.
9 (1972), 253--265.
\bibitem{Cou}
T.J. Courant.
Dirac manifolds.
Trans. Amer. Math. Soc.
319 (1990), no 2, 631--661.
\bibitem{E}
K. Ebrahimi-Fard.
Loday-type algebras and the Rota-Baxter relation.
Lett. Math. Phys.
61 (2002), no 2, 139--147.
\bibitem{FGK1}
K. Ebrahimi-Fard, L. Guo and D. Kreimer.
Integrable renormalization I: the ladder case.
J. Math. Phys.
45 (2004), 3758--3769.
\bibitem{FGK3}
K. Ebrahimi-Fard, L. Guo and D. Kreimer.
Spitzer's Identity and the Algebraic Birkhoff Decomposition in pQFT.
J. Phys. A: Math. Gen.
37 (2004), 11037--11052.
\bibitem{FGK2}
K. Ebrahimi-Fard, L. Guo and D. Kreimer.
Integrable renormalization II: the general case.
Annales Henri Poincare.
6 (2005), 369--395.
\bibitem{Fre}
J.M. Freeman.
On the classification of operator identities.
Studies in Appl. Math.
51 (1972), 73--84.
\bibitem{Kos2}
Y. Kosmann-Schwarzbach.
From Poisson algebras to Gerstenhaber algebras.
Ann. Inst. Fourier (Grenoble).
46 (1996), no 5, 1243--1274.
\bibitem{Ler}
P. Leroux.
Construction of Nijenhuis operators and dendriform trialgebras.
Int. J. Math. Math. Sci.
(2004), no 49-52, 2595--2615.
\bibitem{L4}
J-L. Loday.
Cyclic homology.
Grundlehren der Mathematischen Wissenschaften
(Springer, Berlin).
301 (1992).
\bibitem{L1}
J-L. Loday.
Dialgebras. Dialgebras and related operads.
Lecture Notes in Math (Springer, Berlin).
1763 (2001), 7--66.
\bibitem{Rot1}
G-C. Rota.
Baxter algebras and combinatorial identities. I, II.
Bull. Amer. Math. Soc.
75 (1969), (I) 325--329, (II) 330--334.
\bibitem{Rot2}
G-C. Rota.
Baxter operators, an introduction.
Contemp. Mathematicians (Birkhauser Boston, Boston).
(1995), 504--512.
\bibitem{Se}
P. Severa and A. Weinstein.
Poisson geometry with a 3-form background.
Noncommutative geometry and string theory.
Progr. Theoret. Phys. Suppl.
(2001), no 144, 145--154.
\bibitem{U}
K. Uchino.
Poisson geometry on associative algebras.
In preparation.
\end{thebibliography}
\end{document}